\documentclass[12pt]{amsart}
\advance\textheight\topskip
\allowdisplaybreaks

\usepackage{graphicx}
\usepackage{wrapfig}
\usepackage{mathabx}
\newcommand{\Smash}{\mathbin{\hash}} %%%\#

\newcommand{\leftreg}{\kern1pt{\rightharpoonup}\kern1pt}
\newcommand{\rightreg}{\kern1pt{\leftharpoonup}\kern1pt}

\newcommand{\leftregchi}{\kern1pt{\rightharpoonup}_{\!\!\!\!_{\chi}}\kern1pt}

\usepackage{amsfonts,amssymb}
\usepackage{times}
\usepackage[mathscr]{eucal}
\usepackage{mathptmx}

\usepackage{amsthm}

\newcommand{\Socle}[1]{*{\rule[-1pt]{0pt}{7pt}\bullet}\ar@{}|(.3){#1}[0,0];[0,0]+<0pt,12pt>;}
\newcommand{\Quotient}[1]{*{\rule[-1pt]{0pt}{7pt}\circ}\ar@{}|(.3){#1}[0,0];[0,0]+<0pt,12pt>;}

\newcommand{\Czd}{\oC_{\q}[z,\Dz]}
\newcommand{\Mat}{\mathrm{Mat}}
\renewcommand{\kappa}{\varkappa}
\newcommand{\qfac}[1]{[#1]!\,}
\newcommand{\qint}[1]{[#1]}

\newcommand{\leftii}{\mathbin{\mbox{\small${\vartriangleright}$}}}

\newcommand{\Dp}[1]{\,_{\phantom{h}}^{\underline{#1\kern-.5pt}\kern.5pt}}

\newcommand{\Sinv}{{S^*}^{-1}}

\newcommand{\Dz}{\partial}
\newcommand{\dz}{\underline{dz}}
\newcommand{\dDz}{\underline{d\Dz}}
\newcommand{\dlambda}{\underline{d\lambda}}

\newcommand{\Ri}{R^{(1)}}
\newcommand{\Rii}{R^{(2)}}

\usepackage[curve,matrix,arrow]{xy} %%%,frame
\usepackage[bookmarks=false, hyperfigures=false, a4paper]{hyperref}
\newcommand{\eval}[2]{\langle#1,\,#2\rangle\,}

\newcommand{\DD}{\mathscr{D}}
\newcommand{\HD}{\mathscr{H}}

\newcommand{\bHDB}{\kern4pt\overline{\kern-3pt\mathscr{H}(B)\kern-3pt}\kern4pt}
\newcommand{\bHDBstar}{\kern4pt\overline{\kern-3pt\mathscr{H}(B^*)\kern-3pt}\kern4pt}
\newcommand{\bDDB}{\kern4pt\overline{\kern-3pt\mathscr{D}(B)\kern-3pt}\kern4pt}

\newcommand{\rep}{\mathscr}

\newcommand{\repP}{\rep{P}}

\newcommand{\bref}[1]{\textbf{\ref{#1}}}

\renewcommand{\geq}{\geqslant}
\renewcommand{\leq}{\leqslant}

\newcommand{\tensor}{\otimes}

\newcommand{\q}{\mathfrak{q}}

\newcommand{\UresSL}[1]{\overline{\mathscr{U}}_{\q} s\ell(#1)}
\newcommand{\HresSL}[1]{\overline{\mathscr{H}}_{\q} s\ell(#1)}

\newcommand{\mfrac}[2]{\raisebox{.8pt}{\mbox{\small$\displaystyle\frac{#1}{#2}$}}}
\newcommand{\ffrac}[2]{\raisebox{.5pt}{\mbox{\footnotesize$\displaystyle\frac{#1}{#2}$}}}
\newcommand{\fffrac}[2]{\raisebox{.9pt}{\mbox{\scriptsize$\displaystyle\frac{#1}{#2}$}}}
\newcommand{\half}{%
  \mathchoice{\ffrac{1}{2}}{\frac{1}{2}}{\frac{1}{2}}{\frac{1}{2}}}

\newcommand{\qbin}[2]{\mathchoice%
  {\qbinm{#1}{#2}}{\qbinmm{#1}{#2}}%
  {\qbinmm{#1}{#2}}{\qbinmm{#1}{#2}}}
\newcommand{\qbinm}[2]{\mbox{\footnotesize$\displaystyle
    \genfrac{[}{]}{0pt}{}{#1}{#2}$}}
\newcommand{\qbinmm}[2]{\genfrac{[}{]}{0pt}{}{#1}{#2}}

\newcommand{\dd}{d}
\newcommand{\ddinv}{\dd^{-1}}

\newcommand{\oC}{\mathbb{C}}

\newcommand{\oZ}{\mathbb{Z}}

\numberwithin{equation}{section}
\makeatletter
\@addtoreset{equation}{section}
\@addtoreset{subsubsection}{section}

\def\@secnumfont{\bfseries}
\def\subsubsection{\@startsection{subsubsection}{3}%
  \z@{.5\linespacing\@plus.7\linespacing}{-.5em}%
  {\normalfont\bfseries}}
\def\paragraph{\@startsection{paragraph}{4}%
  \z@\z@{-\fontdimen2\font}%
  \normalfont\bfseries}
\def\subparagraph{\@startsection{subparagraph}{5}%
  \z@\z@{-\fontdimen2\font}%
  \normalfont\bfseries}

\makeatother

\swapnumbers

\theoremstyle{definition}

\begin{document}

\title[HD pair]{%% $\boldsymbol{(\DD(B),\,\HD(B^*))}$
  A Heisenberg double addition to the logarithmic Kazhdan--Lusztig
  duality}

\author[Semikhatov]{A.M.~Semikhatov}%

\address{\mbox{}\kern-\parindent Lebedev Physics Institute
  \hfill\mbox{}\linebreak \texttt{ams@sci.lebedev.ru}}

\begin{abstract}
  For a Hopf algebra $B$, we endow the Heisenberg double $\HD(B^*)$
  with the structure of a module algebra over the Drinfeld double
  $\DD(B)$.  Based on this property, we propose that $\HD(B^*)$ is to
  be the counterpart of the algebra of fields on the quantum-group
  side of the Kazhdan--Lusztig duality between logarithmic conformal
  field theories and quantum groups.  As an example, we work out the
  case where $B$ is the Taft Hopf algebra related to the
  $\UresSL2$ quantum group that is Kazhdan--Lusztig-dual
%%   that underlies the Kazhdan--Lusztig duality
  to $(p,1)$ logarithmic conformal models.  The corresponding pair
  $(\DD(B),\HD(B^*))$ is ``truncated'' to $(\UresSL2,\HresSL2)$, where
%%   $\UresSL2$ is the $2p^3$-dimensional quantum~$s\ell(2)$ and
  $\HresSL2$ is a $\UresSL2$ module algebra that turns out to have the
  form $\HresSL2=\oC_{\q}[z,\Dz]\tensor\oC[\lambda]/(\lambda^{2p}-1)$,
  where $\oC_{\q}[z,\Dz]$ is the $\UresSL2$-module algebra with the
  relations $z^p=0$, $\Dz^p=0$, and $\Dz z = \q-\q^{-1} + \q^{-2}
  z\Dz$.
\end{abstract}

\maketitle

%% \enlargethispage{\baselineskip}

\thispagestyle{empty}
%% \enlargethispage{\baselineskip}

%% \setcounter{tocdepth}{2}%3
%% \vspace*{-24pt}
%% \begin{footnotesize}\addtolength{\baselineskip}{-6pt}
%%   \tableofcontents
%% \end{footnotesize}

\section{Introduction}
The ``logarithmic'' Kazhdan--Lusztig duality\,---\,a remarkable
correspondence between logarithmic conformal field
theories\footnote{It has become impossible to list ``all'' papers on
  logarithmic conformal field theory.  We note the pioneering
  works~\cite{[Gurarie],[RozSal],[K-first],[GK1]}, a prejudiced
  selection~\cite{[FFHST],[FHST],[FGST],[FGST3],[QQ-sl2]}, a
  vertex-operator algebra trend
  in~\cite{[CF],[AM-triplet],[HLZ],[Hu],[AM],[NT],[H]}, and recent
  papers~\cite{[GR2],[FK],[RS],[MR],[SMDR], [PRR],[Ras], [VF],[KR],
    [GRW]}, wherein further references can be found.} and quantum
groups---\,is based on a Drinfeld double construction on the quantum
group side~\cite{[FGST]}.  The starting point is the Hopf algebra $B$
generated by the screening(s) in a logarithmic model and diagonal,
``zero-mode-like'' element(s) (see~\cite{[FGST3],[FGST-q],[Ar-pq]} for
the two-screening case, which is relatively complicated by modern
standards).  The strategy is then to construct the Drinfeld double of
this quantum group and to ``slightly truncate'' it, to produce the
Kazhdan--Lusztig-dual quantum group.  Various aspects of the
``logarithmic'' Kazhdan--Lusztig duality were developed
in~\cite{[FGST2],[FGST3],[FGST-q],[G],[BFGT],[AM],[NT]}.

The resulting correspondence (ranging up to the coincidence) in the
properties of the symmetry algebra of the logarithmic model and the
dual quantum group
%% , such as the representation category and the modular group
%% representation,
is ``circumstantial'' in that it is seen to work nicely in particular
cases, although no general argument for its existence has been
developed or attempted.  That the Drinfeld double of $B$ plays a
crucial role in this correspondence was a serendipitous finding
in~\cite{[FGST]}.  Modulo the ``slight truncation'' mentioned
above,\pagebreak[3]
the Drinfeld double is a counterpart of the symmetry algebra (``the''
triplet~\cite{[K-first],[GK1],[GK3],[FHST],[CF]} or a higher
one~\cite{[FGST3]}) of a given logarithmic conformal field model.

In this paper, we propose another algebraic object that may play a
role in the logarithmic Kazhdan--Lusztig duality, being a counterpart
of the algebra of fields in logarithmic models.  We here mean the
fields describing logarithmic models in manifestly
quantum-group-invariant terms (i.e., ``carrying quantum-group
indices''), as a generalization of the symplectic
fermions~\cite{[K-sy]}.  The necessary algebraic requirement is that
the quantum group act ``covariantly'' on products of fields, which is
expressed as the module algebra axiom
$h\leftii(\varphi\psi)=(h'\leftii\varphi)(h''\leftii\psi)$, where we
use the Sweedler notation $\Delta(h)=h'\tensor h''$ for the coproduct.
We now describe a $\DD(B)$-module algebra that is to play the role of
fields on the algebraic side.

For a Hopf algebra $B$, the Drinfeld double $\DD(B)$ is $B^*\tensor B$
as a vector space.  The same vector space admits another
characteristic algebraic structure, a (semisimple) associative algebra
given by the smash product with respect to the (left) regular action
of $B$ on $B^*$, or, in the established terminology traced back
to~\cite{[AF],[RSts],[Sts]}, a Heisenberg double (see, e.g.,
\cite{[K],[vDvK],[Mil-5]}), specifically, the Heisenberg double
\begin{equation*}
  \HD(B^*)=B^*\Smash B
\end{equation*}
of $B^*$.  The main observation in this paper is that \textit{for any
  Hopf algebra~$B$ with bijective antipode, $\HD(B^*)$ is a
  $\DD(B)$-module algebra}.  This requires introducing a new $\DD(B)$
action (which may be termed ``heterotic'' because it is obtained by
combining, in a sense, a left and a right $\DD(B)$ actions).

As is the case with the Drinfeld double $\DD(B)$, the Heisenberg
double $\HD(B^*)$ turns out to be ``slightly too big'' for the
correspondence with logarithmic models, but for the $2p^3$-dimensional
quantum group $\UresSL2$ at the $2p$th root of unity dual to the
$(p,1)$ logarithmic conformal models, the corresponding $\HD(B^*)$
nicely allows a ``truncation'' to a $2p^3$-dimensional
$\UresSL2$-module algebra.

We prove the general statement in Sec.~\ref{sec:thm} and detail the
$\UresSL2$ example in Sec.~\ref{sec:sl2}.  The definition of the
Drinfeld double is recalled in Appendix~\ref{app:D-double}.  In
Appendix~\ref{app:lcftrem}, we collect some motivation coming from
logarithmic conformal field theories.

\section{$\HD(B^*)$ as a $\DD(B)$-module algebra}\label{sec:thm}
Let $B$ be a Hopf algebra.  In this section, we make $\HD(B^*)$ into a
$\DD(B)$-module algebra.  For this, we combine two well-known $\DD(B)$
actions, which can be taken from different sources, among which we
prefer the beautiful paper~\cite{[P]}.

\subsection{}
We use the ``tickling'' notation for the left and right regular
actions: for a Hopf algebra $H$, its left and right regular actions on
$H^*$ are respectively given by $h\leftreg
\beta=\beta(?h)=\eval{\beta''}{h}\beta'$ and $\beta\rightreg
h=\beta(h?)$, where $\beta\in H^*$ and $h\in H$. It follows that $H^*$
is an $H$-bimodule under these actions (and $\eval{~}{~}$ is the
evaluation).  We also have the left and right actions of $H^*$ on
$H$, $\beta\leftreg a=\eval{\beta}{a''}a'$ and
$a\rightreg\beta=\eval{\beta}{a'}a''$.
%% We use $\eval{\beta}{a}$ and $\beta(a)$ as synonyms.

\subsection{}We recall that the Heisenberg double $\HD(B^*)$ is the
smash product $B^*\Smash B$ with respect to the left regular action of
$B$ on $B^*$, which means that the composition in $\HD(B^*)$ is
given~by
\begin{equation}\label{H-comp}
  (\alpha\Smash a)(\beta\Smash b)=
  \alpha(a'\leftreg\beta)\Smash a'' b,
  \qquad
  \alpha,\beta\in B^*,\quad
  a,b\in B.
\end{equation}
%% As an aside, we note a property of the Heisenberg double known from
%% \cite{[Lu-alg]}: $B^*\Smash B$ is a Hopf algebroid over~$B^*$.

We now describe the $\DD(B)$ action on $\HD(B^*)$ making it into a
$\DD(B)$-module algebra.

First, the $\DD(B)$ action on $B^*$\,---\,the first factor in
$\HD(B^*)=B^*\Smash B$\,---\,is given by the restriction of the left
regular action of $\DD(B)$ on $\DD(B)^*\cong B\tensor B^*$, which
is~\cite{[Lu-double]}
\begin{align*}
  (\mu\tensor m)\leftreg (a\tensor\alpha)
%%   &= \eval{\mu''}{a''} \eval{\alpha''}{n}\,
%%   a' \tensor \mu'''\alpha' \Sinv(\mu')
%%   \\
  &= (\mu''\leftreg a)\tensor
  \mu'''(m\leftreg\alpha)\Sinv(\mu').
\end{align*}
Restricting this to $1\tensor B^*$ gives
\begin{equation}\label{action-i}
  (\mu\tensor m)\leftreg\alpha=\mu''(m\leftreg\alpha)\Sinv(\mu'),
  \qquad
  \mu\tensor m\in\DD(B),\quad
  \alpha\in B^*,
\end{equation}
under which $B^*$ is an $R$-commutative $\DD(B)$-module
algebra~\cite{[Lu-alg]} (also see~\cite{[P]}).\footnote{An algebra $A$
  carrying an action of a quasitriangular Hopf algebra $H$ is called
  $R$-commutative, or quantum commutative~\cite{[CW],[P]} if $a b =
  (\Rii.b)(\Ri.a)$ for all $a,b\in A$, where the dot denotes the
  action and $R=\Ri\tensor\Rii\in H\tensor H$ is the universal
  $R$-matrix.}

Second, the $\DD(B)$ action on $B$ is obtained by restricting the
right regular action of $\DD(B)$ on $\DD(B)^*\cong B\tensor B^*$ to
$B\tensor\varepsilon$ and using the antipode to convert it into a left
action~\cite{[Z]}.  With the right regular action of $\DD(B)$ on
$\DD(B)^*$ given by~\cite{[Lu-double],[P]}
\begin{equation*}
  (a\tensor\alpha)\rightreg(\mu\tensor m)
%%   = \mu(a')\alpha'(m'')
%%   S^{-1}(m''') a'' m'\tensor\alpha''
  =S^{-1}(m''')(a\rightreg\mu)m'\tensor(\alpha\rightreg m''),
\end{equation*}
its restriction to $B$ is $a\rightreg(\mu\tensor m)
=S^{-1}(m'')(a\rightreg\mu)m'$.  Replacing $\mu\tensor m$ here with
$S_{_{\DD}}(\mu\tensor m) =(S(m''')\leftreg\Sinv(\mu)\rightreg
m')\tensor S(m''),$ we readily calculate $a\rightreg
S_{_{\DD}}(\mu\tensor m)= \eval{\Sinv(\mu)}{m'a'S(m'''')}m'' a''
S(m''')$, which defines the left action~\cite{[Z]}
\begin{equation}\label{action-ii}
  (\mu\tensor m)\leftii a
  =(m' a S(m''))\rightreg\Sinv(\mu),
  \qquad
  \mu\tensor m\in\DD(B),\quad a\in B,
\end{equation}
under which $B$ is an $R$-commutative $\DD(B)$-module algebra (also
see~\cite{[P]}).

We now define a $\DD(B)$ action on $\HD(B^*)$, also denoted by
${}\leftii$, simply by setting\footnote{The coproduct
  in~\eqref{action-full} refers to $\DD(B)$, and hence, in accordance
  with the Drinfeld double construction, $(\mu\tensor
  m)'\tensor(\mu\tensor m)'' =(\mu''\tensor m')\tensor(\mu'\tensor
  m'')$, with the coproducts of $B^*$ and $B$ in the right-hand side.}
\begin{equation}\label{action-full}
  (\mu\tensor m)\leftii (\alpha\Smash a)
  =  \bigl((\mu\tensor m)'\leftreg\alpha\bigr)\Smash
  \bigl((\mu\tensor m)''\leftii a\bigr),
\end{equation}
that is,
\begin{align*}
  (\mu\tensor m)\leftii(\alpha\Smash a)
  =\mu'''(m'\leftreg\alpha)\Sinv(\mu'')
  \Smash\bigl((m'' a S(m'''))\rightreg\Sinv(\mu')\bigr),\\[-2pt]
  \notag
  \mu\tensor m\in\DD(B),\quad \alpha\Smash a\in\HD(B^*),
\end{align*}
and prove that $\HD(B^*)$ is then a $\DD(B)$-module algebra.  Because
each factor in $\HD(B^*)= B^*\Smash B$ is already a $\DD(B)$-module
algebra, it suffices to show that
\begin{equation*}
  \bigl(\!(\mu\tensor m)'\leftii(\varepsilon\Smash a)\!\bigr)
  \bigl(\!(\mu\tensor m)''\leftii(\beta\Smash 1)\!\bigr)
  =(\mu\tensor m)\leftii
  \bigl((a'\leftreg\beta)\Smash a''\bigr).
\end{equation*}
We evaluate the left-hand side:
\begin{multline*}
  \bigl(\!(\mu''\tensor m')\leftii(\varepsilon\Smash a)\!\bigr)
  \bigl(\!(\mu'\tensor m'')\leftii(\beta\Smash 1)\!\bigr)
  \\
  \begin{aligned}
    &=\bigl((\mu''\tensor m')\leftii a\bigr)'
    \leftreg
    (\mu'\tensor m''\leftreg\beta)
    \Smash
    \bigl((\mu''\tensor m')\leftii a\bigr)''
    \\
    &=\bigl((m^{(1)} a' S(m^{(4)}))\rightreg\Sinv(\mu'')\bigr)
    \leftreg(\mu'\tensor m^{(5)}\leftreg\beta)
    \Smash
    m^{(2)} a'' S(m^{(3)})
    \\*[-2pt]
    &\phantom{{}={}}
    \text{\mbox{}\hfill(because $\Delta((\mu\tensor m)\leftii a)
      =\bigl(m' a' S(m'''')\rightreg\Sinv(\mu)\bigr)
      \tensor m'' a'' S(m''')$)}
    \\[2pt]
    &=(\mu'\tensor m^{(5)}\leftreg\beta)'\Smash m^{(2)} a'' S(m^{(3)})
    \eval{\Sinv(\mu'')(\mu'\tensor m^{(5)}\leftreg\beta)''}{
      m^{(1)} a' S(m^{(4)})}\\*[-2pt]
    &\phantom{{}={}}
    \text{(simply because $(a\rightreg\alpha)\leftreg\beta
      =\beta'\eval{\alpha\beta''}{a}$)}
    \\[2pt]
    &=\mu^{(3)}\beta'\Sinv(\mu^{(2)})\Smash m^{(2)} a'' S(m^{(3)})
    \\*
    &\qquad\qquad\qquad{}\times
    \eval{\Sinv(\mu^{(5)})\mu^{(4)}(m^{(5)}\leftreg\beta'')
    \Sinv(\mu^{(1)})}{m^{(1)} a' S(m^{(4)})}
  \\*[-2pt]
  &\phantom{{}={}}
  \text{(because $\Delta((\mu\tensor m)\leftreg\beta)
    =\mu'''\beta'\Sinv(\mu'')\tensor\mu''''(m\leftreg\beta'')
    \Sinv(\mu')$)}
  \\[2pt]
  &=\mu^{(3)}\beta'\Sinv(\mu^{(2)})\Smash m^{(3)}a''' S(m^{(4)})
  \\*
  &\qquad\qquad\qquad{}\times
  \eval{m^{(7)}\leftreg\beta''}{m^{(1)} a' S(m^{(6)})}
  \eval{\Sinv(\mu^{(1)})}{m^{(2)} a'' S(m^{(5)})}
  \\
  &=\mu'''(m'a'\leftreg\beta)\Sinv(\mu'')
  \Smash
  \bigl((m'' a'' S(m'''))\rightreg\Sinv(\mu')
  \bigr)
  \end{aligned}
  \\=\bigl((m'\tensor\mu'')\leftreg(a'\leftreg\beta)\bigr)
  \Smash
  \bigl((m''\tensor\mu')\leftii a''\bigr),
\end{multline*}
which is the desired result.

\subsection{Remark} As already noted, each of the
subalgebras $B^*\tensor 1$ and $\varepsilon\tensor B$ in $\HD(B^*)$ is
known to be $R$-commutative with respect to the corresponding
action~\eqref{action-i} or~\eqref{action-ii} of $\DD(B)$.  But
$\HD(B^*)$ is \textit{not} $R$-commutative with respect to the action
in~\eqref{action-full} in general: the $R$-commutativity axiom is
satisfied for only ``half'' the cross-relations,
\begin{equation*}%%%\label{cross-qc}
  \bigl(\Rii\leftii(\varepsilon\Smash b)\bigr)
  \bigl(\Ri\leftii(\alpha\Smash 1)\bigr)
  =(\alpha\Smash 1) (\varepsilon\Smash b)=\alpha\Smash b,
\end{equation*}
but not for the other half: $\bigl(\Rii\leftii(\beta\Smash 1)\bigr)
\bigl(\Ri\leftii(\varepsilon\Smash a)\bigr)\neq (\varepsilon\Smash
a)(\beta\Smash 1)$ in general.  
%% For completeness, we now show~\eqref{cross-qc}, by evaluating the
%% left-hand side:
%% \begin{align*}
%%   \bigl(\varepsilon\Smash(e^I\leftii b)\bigr)
%%   \bigl((e_I\leftreg\alpha)\Smash 1\bigr)
%%   &=\bigl(\varepsilon\Smash(b\rightreg\Sinv(e^I))\bigr)
%%   \bigl((e_I\leftreg\alpha)\Smash 1\bigr)
%%   \\[-2pt]
%%   &=\bigl((b\rightreg\Sinv(e^I))'\leftreg
%%   e_I\leftreg\alpha\bigr)\Smash(b\rightreg\Sinv(e^I))''
%%   \\[-2pt]
%%   &=\bigl((b'\rightreg\Sinv(e^I))\leftreg
%%   e_I\leftreg\alpha\bigr)\Smash b''
%%   \\[-2pt]
%%   &=\eval{\Sinv(e^I)(e_I\leftreg\alpha)''}{b'}
%%   (e_I\leftreg\alpha)' \Smash b''
%%   \\[-2pt]
%%   &=\eval{\Sinv(e^I)(e_I\leftreg\alpha'')}{b'} \alpha' \Smash b''
%%   \\[-2pt]
%%   &=\eval{\Sinv(e^I)}{b'}\eval{e_I\leftreg\alpha''}{b''} \alpha'
%%   \Smash b''
%%   \\[-2pt]
%%   &=\eval{e^I}{S^{-1}(b')}\eval{\alpha''}{b''e_I} \alpha'
%%   \Smash b''
%%   \\[-2pt]
%%   &=\eval{\alpha''}{b''S^{-1}(b')} \alpha' \Smash b''=\alpha\Smash b.
%% \end{align*}

\section{The $(\UresSL2,\,\HresSL2)$ pair}\label{sec:sl2}
In this section, we consider the pair $(\DD(B),\HD(B^*))$ for the Taft
Hopf algebra $B$ that underlies the Kazhdan--Lusztig correspondence
with the $(p,1)$ logarithmic conformal field theory models.  By
``truncation,'' $\DD(B)$ yields the $\UresSL2$ quantum group that is
Kazhdan--Lusztig-dual to the $(p,1)$ logarithmic models.  This quantum
group first appeared in~\cite{[AGL]} and was rediscovered, together
with its role in the Kazhdan--Lusztig correspondence,
in~\cite{[FGST]}; its further properties were considered
in~\cite{[FGST2],[S-q],[Ar],[FHT],[S-U]} and, notably, recently
in~\cite{[KS]} (also see~\cite{[Erd]} for a somewhat larger quantum
group).  We recall this in~\bref{the-B}.  We evaluate $\HD(B^*)$
in~\bref{HD-relations}, and in~\bref{to-sl2} ``truncate''
$(\DD(B),\HD(B^*))$ to a pair $(\UresSL2,\,\HresSL2)$, where
$\HresSL2$ is a $\UresSL2$-module algebra.  Its structure is detailed
in~\bref{sec:structure}.

\subsection{$\DD(B)$ for the $4p^2$-dimensional Taft Hopf algebra
  $B$}\label{the-B}
For an integer $p\geq2$, we set
\begin{gather}\label{the-q}
  \q=e^{\frac{i\pi}{p}}
\end{gather}
and recall some of the results in~\cite{[FGST]}.

\subsubsection{The Taft Hopf algebra $B$}
Let 
\begin{equation*}
  B=\mathrm{Span}(E^m k^n),\quad
  0\leq m\leq p-1,\quad 0\leq n\leq 4p-1,
%%   E^m k^n &= E^m k^n,
\end{equation*}
be the $4p^2$-dimensional Hopf algebra generated by~$E$ and~$k$ with
the relations
\begin{gather}\label{prod-B}
  k E =\q E k,\quad E^p=0,\quad k^{4p}=1,
\end{gather}
and with the comultiplication, counit, and antipode given by
\begin{gather}\label{coalgebra-B}
  \begin{gathered}
    \Delta(E)=1\otimes E+ E\otimes k^2,\quad
    \Delta(k)= k\otimes k,\\
    \epsilon(E)=0,\quad\epsilon(k)=1,\\
    S(E)=- E k^{-2},\quad S(k)= k^{-1}.
  \end{gathered}
\end{gather}

\subsubsection{$B^*$ and $\DD(B)$}\label{sec:B-dual}
We next introduce elements $F,\varkappa\in B^*$ as
\begin{equation*}
  \eval{F}{E^m k^n}=\delta_{m,1}\ffrac{\q^{-n}}{\q-\q^{-1}},
  \qquad
  \eval{\varkappa}{E^m k^n}=\delta_{m,0}\q^{-n/2}.
\end{equation*}
Then~\cite{[FGST]}
\begin{align*}
  B^*&=\mathrm{Span}(F^a\varkappa^b),\quad
  0\leq a\leq p-1,\quad 0\leq b\leq 4p-1.
\end{align*}
Moreover, straightforward calculation shows~\cite{[FGST]} that the
Drinfeld double $\DD(B)$ (see Appendix~\ref{app:D-double}) is the Hopf
algebra generated by $E$, $F$, $k$, and $\varkappa$ with the relations
given~by
\begin{itemize}
\item[i)] relations~\eqref{prod-B} in~$B$,
\item[ii)] the relations
  \begin{equation*}%%%\label{prod-Bstar}
    \varkappa F=\q F\varkappa,\quad F^p=0,\quad
    \varkappa^{4p}=1
  \end{equation*}
  in $B^*$, and
\item[iii)] the cross-relations
  \begin{gather}\label{cross}
    k\varkappa=\varkappa k,\quad k F k^{-1}=\q^{-1} F,\quad
    \varkappa E\varkappa^{-1}=\q^{-1} E,\quad
    [E, F]=\mfrac{ k^2-\varkappa^2}{\q-\q^{-1}}.
  \end{gather}
\end{itemize}
Here, in accordance with writing $\DD(B)=B^*\tensor B$, $E$ and $k$
are of course understood as $\varepsilon\tensor E$ and
$\varepsilon\tensor k$, and $F$ and $\varkappa$ as $F\tensor 1$ and
$\varkappa\tensor 1$.  Then, for example, the last relation
in~\eqref{cross} is to be rewritten as\pagebreak[3]
\begin{equation*}
  (\varepsilon\tensor E)(F\tensor 1) =F\tensor E +
  \ffrac{1}{\q-\q^{-1}}
  \varepsilon\tensor k^2 - \ffrac{1}{\q-\q^{-1}}\varkappa^2\tensor 1.
\end{equation*}
Dropping the ${}\tensor{}$ inside $\DD(B)$ again, we have the
Hopf-algebra structure
$(\Delta_{_{\DD}},\varepsilon_{_{\DD}},S_{_{\DD}})$ given
by~\eqref{coalgebra-B} and
\begin{gather*}
%%   \Delta(E)=1\tensor E+ E\tensor k^2,\quad
%%     \Delta(k)= k\tensor k,\quad
%%     \epsilon(E)=0,\quad
%%     \epsilon(k)=1,\\
  \Delta_{_{\DD}}(F)=\varkappa^2\tensor F+ F\tensor1,\quad
  \Delta_{_{\DD}}(\varkappa)=\varkappa\tensor\varkappa,\quad
  \epsilon_{_{\DD}}(F)=0,\quad
  \epsilon_{_{\DD}}(\varkappa)=1,
%%     \\
%%     S_{_{\DD}}(E)=- E k^{-2},\quad
%%     S_{_{\DD}}(k)= k^{-1},
  \\
  S_{_{\DD}}(F)=-\varkappa^{-2} F,
  \quad S_{_{\DD}}(\varkappa)=\varkappa^{-1}
\end{gather*}
(we reiterate that the coalgebra structure on $\DD(B)$ is the direct
product of those on $B^{*\mathrm{cop}}$ and $B$).  
%% It also follows that
%% \begin{align*}
%%   \Delta(E^{m}) &= \sum_{s=0}^{m}\q^{-s(s - m)}\qbin{m}{s} E^{s}\tensor
%%   E^{m - s}k^{2s},
%%   \\
%%   \Delta_{_{\DD}}
%%   (F^{m})
%%   &= \sum_{s=0}^{m} \q^{-s(s - m)}\qbin{m}{s}
%%   F^{m - s}
%%   \varkappa^{2s} %%% k^{-2s}
%%   \tensor F^{s}.
%% \end{align*}
%% Some other formulas pertaining to the explicit construction of
%% $\DD(B)$ are given in~\bref{double-specific}.

\subsection{The Heisenberg double $\HD(B^*)$}\label{HD-relations}
For the above $B$, \ $\HD(B^*)$ is spanned by
\begin{equation}\label{HBstar-span}
  F^a\varkappa^b\Smash E^c k^d,\qquad a,c=0,\dots,p-1,\quad
  b,d\in\oZ/(4p\oZ),
%%   =0,\dots,4p-1
\end{equation}
%% (with the range of $b$ and $d$ understood as $\oZ/(4p\oZ)$), 
where $\varkappa^{4p}=1$, $k^{4p}=1$, $F^p=0$, and $E^p=0$.

\subsubsection{The composition law}
To evaluate the product in $\HD(B^*)$, defined in~\eqref{H-comp}, we
first write the left regular action of $B$ on $B^*$, $b\leftreg\beta
%%%=\beta'\eval{\beta''}{b} 
=\beta''_{_{\DD}} \eval{\beta'_{_{\DD}}}{b}$:
\begin{gather}\label{D-action:Ek}
  E^m k^n\leftreg (F^a\varkappa^b)
  =\qbin{a}{m}\ffrac{\qfac{m}}{(\q - \q^{-1})^m}\,
  \q^{-(b+2a)\frac{n}{2} - m(a + b) + \half m(m+1)}
  F^{a-m}\varkappa^b.
\end{gather}
%% In particular,
%% \begin{align*}
%%   E \leftreg (F^a\varkappa^b)
%%   &=\ffrac{\qint{a}}{\q - \q^{-1}}\,
%%   \q^{1-a-b}
%%   F^{a-1}\varkappa^b,\\%%%[-6pt]
%% %%%   \intertext{and}
%%   k\leftreg (F^a\varkappa^b)
%%   &= \q^{-a - \frac{b}{2}}
%%   F^{a}\varkappa^b.
%% \end{align*}
It then follows that
\begin{multline}\label{prod-in-smash}
  (\varepsilon\Smash E^m k^n)(F^a\varkappa^b\Smash 1)
%%   =(E^m k^n)'\leftreg(F^a\varkappa^b)\Smash (E^m k^n)''
  \\
  {}=
  \sum_{s\geq 0}\q^{-\half s(s-1)}\qbin{m}{s}
  \qbin{a}{s}\ffrac{\qfac{s}}{(\q - \q^{-1})^s}\,
  \q^{-(b+2a)\frac{n}{2} + s(m - a - b)}
  F^{a-s}\varkappa^b
  \Smash E^{m - s}k^{2s+n}
\end{multline}
(the sum is limited above by $\min(m,a)$ due to the binomial
coefficient vanishing).  In particular,
\begin{gather*}
  (\varepsilon\Smash E k^n)(F\varkappa^b\Smash 1)
%%   =
%%   \sum_{s=0}^{1}\ffrac{1}{(\q - \q^{-1})^s}\,
%%   \q^{-(b + 2)\frac{n}{2} + s(-b)}
%%   F^{1-s}\varkappa^b
%%   \Smash E^{1 - s}k^{2s+n}\\
  =
  \q^{-(b + 2)\frac{n}{2}}
  F\varkappa^b
  \Smash Ek^{n}
  +\ffrac{1}{\q - \q^{-1}}\,
  \q^{-(b + 2)\frac{n}{2} - b}
  \varkappa^b \Smash k^{n+2},
\end{gather*}
and also $(\varepsilon\Smash k)(\varkappa\Smash 1)
=\q^{-\half}\varkappa\Smash k$, \;$(\varepsilon\Smash k)(F\Smash 1)
=\q^{-1}F\Smash k$, and \;$(\varepsilon\Smash E)(\varkappa\Smash
1)={}$\linebreak[0]$\varkappa\Smash E$.  For the future reference, we
write the general case, obtained from~\eqref{prod-in-smash}
immediately:
\begin{multline}\label{prod-in-smash-full}
  (F^r \varkappa^s\Smash E^m k^n)(F^a\varkappa^b\Smash E^c k^d)
  \\
  {}= \sum_{u\geq 0}\q^{-\half u(u-1)}\qbin{m}{u}
  \qbin{a}{u}\ffrac{\qfac{u}}{(\q - \q^{-1})^u}\,
%%   \q^{-(b+2a)\frac{n}{2} + u(m - a - b) + s(a-u) + (2u+n)c}
  \q^{-\half b n + c n + a (s - n) +  u(2 c - a - b  + m - s)}
  \\*[-6pt]
  {}\times F^{a+r-u}\varkappa^{b+s} \Smash E^{m + c - u}k^{n+d+2u}.
\end{multline}
(This is an associative product for generic $q$ as well.)

\subsubsection{The $\DD(B)$ action}\label{DD-action}We next evaluate
the $\DD(B)$ action on $\HD(B^*)$.

The $\DD(B)$ action on $B^*$ in~\eqref{action-i}, rewritten in terms
of the comultiplication and antipode of the double,
\begin{equation*}
  (\mu\tensor m)\leftreg\alpha=
  \eval{\alpha'_{_{\DD}}}{m}\mu'_{_{\DD}}\alpha''_{_{\DD}}
  S_{_{\DD}}(\mu''_{_{\DD}}),
  \qquad
  \mu=F^i\varkappa^j,\quad
  m=E^m k^n,
\end{equation*}
factors into the action of $\varepsilon\tensor m$
in~\eqref{D-action:Ek} times the action of $\mu\tensor 1$ given by
\begin{equation*}%%%\label{D-action:Fkappa}
  F^i\varkappa^j\leftreg (F^a\varkappa^b)
  =\q^{\frac{i}{2}(i - 1 + b) + a(i + j)}
  (-1)^{i}
  (\q-\q^{-1})^{i}
  \prod_{\ell=1}^{i}\,\qint{\ell + a - 1 + \fffrac{b}{2}}\;
  F^{i+a}\varkappa^{b}.
\end{equation*}
%% In particular,
%% \begin{align*}
%%   F\leftreg(F^a\varkappa^b)
%%   &=(1-\q^{2 a + b}) F^{a+1}\varkappa^b
%%   \\
%% %%   \intertext{and}
%%   \varkappa\leftreg(F^a\varkappa^b)
%%   &=\q^{a}
%%   F^{a}\varkappa^b.
%% \end{align*}

The $\DD(B)$ action on $B$ in~\eqref{action-ii}, $(\mu\tensor
m)\leftii a =(m' a S(m''))\rightreg S_{_{\DD}}(\mu)$, with
$\mu=F^i\varkappa^j$ and $m=E^m k^n$, factors through the adjoint
action of $\varepsilon\tensor m\in\varepsilon\tensor B$,
\begin{equation*}
  E^m k^n\leftii(E^a k^b) = \q^{a n + \half m(1 - m + b)}
  (\q-\q^{-1})^m
  \Bigl(\prod_{\ell=1}^m\qint{\ell-1-\fffrac{b}{2}}\Bigr)
  E^{a+m} k^{b-2m},
\end{equation*}
and the action of $\mu\tensor 1\in B^*\tensor 1$, given by
$\mu\leftii{}a
%%% =a\rightreg\Sinv(\mu)=\eval{\Sinv(\mu)}{a'}a''
=\eval{S_{_{\DD}}(\mu)}{a'}a''$:
\begin{equation*}
  F^i\varkappa^j\leftii(E^a k^b)=
  (-1)^i\qbin{a}{i}\ffrac{\qfac{i}}{(\q-\q^{-1})^i}\,
  \q^{\frac{b j}{2} - \half i(i+1) +i(j+a)}
  E^{a-i} k^{2i+b}.
\end{equation*}

The action in~\eqref{action-full} is therefore given by
\begin{align*}
  E^m\leftii(F^a\varkappa^b\Smash E^c k^d)
  &=\q^{-\half m(m-1)}\sum_{s\geq 0}\q^{-s^2 + 2 s m + s(2c - a - b) +\half d(m-s)}
  \qbin{m}{s}\qbin{a}{s}\qfac{s}\\[-4pt]
  &\quad{}\times
  \Bigl(\prod_{\ell=1}^{m-s}\qint{\ell-1-\fffrac{d}{2}}\Bigr)
  (\q-\q^{-1})^{m-2s}
  F^{a-s}\varkappa^b\Smash E^{c+m-s} k^{d - 2 m + 2 s},
  \\
  k\leftii(F^a\varkappa^b\Smash E^c k^d)
  &=\q^{-a + c - \frac{b}{2}}(F^a\varkappa^b\Smash E^c k^d),
  \\
  \varkappa\leftii(F^a\varkappa^b\Smash E^c k^d)
  &=\q^{a + \frac{d}{2}}(F^a\varkappa^b\Smash E^c k^d),
  \\
  F^i\leftii(F^a\varkappa^b\Smash E^c k^d)
  &=\q^{\half i(i-1)}\sum_{s\geq 0}(-1)^i\q^{-s^2}\qbin{i}{s}\qbin{c}{s}\qfac{s}
  \q^{\half b(i-s) + a i + a s + s c}\\[-2pt]
  &\quad{}\times
  \Bigl(\prod_{\ell=1}^{i-s}\qint{\ell + a - 1 + \fffrac{b}{2}}\Bigr)
  (\q - \q^{-1})^{i - 2 s}
  F^{a  + i - s}\varkappa^b\Smash E^{c-s} k^{d + 2 s}.
\end{align*}

\subsection{From $\DD(B)$ to $\UresSL2$}\label{to-sl2}
The ``truncation'' whereby $\DD(B)$ yields $\UresSL2$~\cite{[FGST]}
consists of two steps: first, taking the quotient
\begin{gather}\label{quotient}
  \bDDB=\DD(B)/(\varkappa k - 1)
\end{gather}
by the Hopf ideal generated by the central element $\varkappa\tensor k
- \varepsilon\tensor 1$ and, second, identifying $\UresSL2$ as the
subalgebra in $\bDDB$ spanned by $F^{\ell} E^{m} k^{2n}$ (tensor
product omitted) with $\ell,m=0,\dots,p-1$ and $n=0,\dots,2p-1$.  It
follows from the above formulas for $\Delta$ and from formulas for the
antipode that $\UresSL2$ is a Hopf algebra.\footnote{It is actually a
  ribbon and (slightly stretching the definition) factorizable Hopf
  algebra~\cite{[FGST],[FGST2],[S-q]}\,---\,the properties playing a
  crucial role in the Kazhdan--Lusztig correspondence.}

In $\HD(B^*)$, dually, we take a subalgebra and then a quotient, as
follows.

First, dually to taking the quotient in~\eqref{quotient}, we identify
the subspace $\bHDBstar\subset\HD(B^*)$ on which $\varkappa\tensor
k\in\DD(B)$ acts by unity.  It follows from the above formulas for the
$\DD(B)$ action that
\begin{align*}
  \bHDBstar&=\mathrm{Span}(\Psi^{a,b,c}),
  \qquad a,c=0,\dots,p-1,\quad
  b\in\oZ/(4p\oZ),
  \\
  \Psi^{a,b,c}&=F^a \varkappa^b\Smash E^c k^{b-2c}.
\end{align*}
Two nice properties immediately follow:
from~\eqref{prod-in-smash-full}, $\bHDBstar$ is a subalgebra, and
from~\bref{DD-action}, the $\DD(B)$ action restricts to $\bHDBstar$.
%% \begin{small}
%%   Explicitly,
%%   \begin{align*}
%%     E\leftii(F^a\varkappa^b\Smash E^c k^d)
%%     &=(1-\q^{d})F^a\varkappa^b\Smash E^{c+1} k^{d-2}
%%     +\ffrac{\qint{a}}{\q-\q^{-1}}\q^{1-a-b+2c}
%%     F^{a-1}\varkappa^b\Smash E^c k^d,
%%     \\
%%     F\leftii(F^a\varkappa^b\Smash E^c k^d)
%%     &=(1-\q^{2a+b})F^{a+1}\varkappa^b\Smash E^c k^d
%%     -\q^{2a+c-1}\ffrac{\qint{c}}{\q-\q^{-1}}
%%     F^a\varkappa^b\Smash E^{c-1}k^{d+2}.
%%   \end{align*}%
%% \end{small}%

Second, dually to the restriction $\UresSL2\subset\bDDB$, we take a
quotient of $\bHDBstar$.  It follows from
$k^2\leftii(F^a\varkappa^b\Smash E^c k^d)= \q^{-2a - b +
  2c}\,F^a\varkappa^b\Smash E^c k^d$ that the eigenvalues of $(k^2)^b$
are not all different for $b\in\oZ/(4p\oZ)$; we can impose the
additional relation $\varkappa^{2p}\Smash k^{2p}=1$
in~$\bHDBstar$,\footnote{The element $\Lambda=\varkappa^{2p}\Smash
  k^{2p}$ is central in $\bHDBstar$, which suffices for our purposes,
  although it is not central in $\HD(B^*)$, where
  $\Lambda\,F^a\varkappa^b\Smash E^c k^d =(-1)^b
  F^a\varkappa^{b+2p}\Smash E^c k^{d+2p}$ and $F^a\varkappa^b\Smash
  E^c k^d\,\Lambda =(-1)^d F^a\varkappa^{b+2p}\Smash E^c k^{d+2p}$.}
i.e., pass to the quotient by the relations
\begin{equation*}
  \Psi^{a,b+2p,c}=(-1)^b\,\Psi^{a,b,c}.
\end{equation*}
This defines the $2p^3$-dimensional algebra $\HresSL2$, which is a
$\UresSL2$ module algebra.

\subsection{The structure of $\HresSL2$}\label{sec:structure}

\subsubsection{}\label{change-basis}Being a semisimple associative
algebra, a Heisenberg double decomposes into matrix algebras.  For our
$\HD(B^*)$, we choose the generators as $(\varkappa, z, \lambda,
\Dz)$, where $\varkappa$ is understood as $\varkappa\Smash1$ and we
set
\begin{align*}
  z &= -(\q-\q^{-1}) \varepsilon\Smash E k^{-2},\\
  \lambda &=\varkappa \Smash k,\\
  \Dz &= (\q-\q^{-1}) F\Smash 1.
\end{align*} 
The relations in $\HD(B^*)$ are then equivalent to
\begin{gather}
  \varkappa^{4p}=1, \qquad\lambda^{4p}=1,\\
  z^p=0,\qquad \Dz^p=0,\label{zp-and-dp}
  \\
  \Dz z = (\q-\q^{-1}) 1 +  \q^{-2} z\Dz,\label{dz-to-zd}\\
    \lambda z = z \lambda,\qquad   \lambda \Dz = \Dz \lambda,\\
  \varkappa z = \q^{-1} z \varkappa,\quad
  \varkappa \lambda = \q^{\half} \lambda \varkappa,\quad
  \varkappa \Dz = \q \Dz \varkappa
\end{gather}
(where the unity in~\eqref{dz-to-zd} is of course $\varepsilon\Smash
1$ in the detailed nomenclature used above).  Clearly, $\lambda$, $z$,
and $\Dz$ generate a subalgebra, which is in fact $\bHDBstar$.  Its
quotient by $\lambda^{2p}=1$ gives $\HresSL2$.  It follows that as an
associative algebra,\pagebreak[3]
\begin{equation*}
  \HresSL2 = \Czd\tensor(\oC[\lambda]/(\lambda^{2p}-1)),
\end{equation*}
where $\Czd$ is the $p^2$-dimensional algebra defined by
relations~\eqref{zp-and-dp} and~\eqref{dz-to-zd}.  It is indeed
isomorphic to the full matrix algebra $\Mat_p(\oC)$~\cite{[S-U]} (also
see~\cite{[FIK]}); hence, $\HresSL2
\cong\Mat_p\bigl(\oC[\lambda]/(\lambda^{2p}-1)\bigr)$.

The $\UresSL2$ action on the new generators of $\HD(B^*)$ is readily
seen to be given by
%% \footnote{Also, $E\leftii k=\fffrac{\q}{\q+1}z\,k$, \ $k^2\leftii
%%   k= k$, \ $F\leftii k=0$.}
\begin{alignat*}{3}
  E\leftii\varkappa &= 0,
  &
  k^2\leftii\varkappa &= \q^{-1}\varkappa,
  &
  F\leftii\varkappa&=-\ffrac{\q}{\q+1}\,\Dz\varkappa,
  \\
%%   E\leftii\lambda&=\ffrac{1}{\q+1}\,\lambda\,z,
  E\leftii\lambda^n&=\q^{-\frac{n}{2}}\qint{\fffrac{n}{2}}\,\lambda^n\,z,
  &
%%   k^2\leftii\lambda&=\q^{-1}\lambda,
  k^2\leftii\lambda^n&=\q^{-n}\lambda,
  &
%%   F\leftii\lambda&=-\ffrac{\q}{\q+1}\Dz\,\lambda,\quad
  F\leftii\lambda^n&=-\q^{\frac{n}{2}}\qint{\fffrac{n}{2}}\,\lambda^n\,\Dz,
  \\
  E\leftii z^m &=-\q^m \qint{m} z^{m+1},
  &
  k^2\leftii z^m &= \q^{2m}\,z^m,
  &
  F\leftii z^m &= \qint{m} \q^{1-m}\,z^{m-1},
  \\
  E\leftii \Dz^n &= \q^{1-n}\qint{n}\Dz^{n-1},\quad
  &
  k^2\leftii \Dz^n &= \q^{-2 n} \Dz^n,\quad
  &
  F\leftii \Dz^n &= -\q^n \qint{n} \Dz^{n+1}
\end{alignat*}
%% (the action on $\varkappa$ and $\Dz$ reduces to the ${}\leftreg{}$
%% above, but we use ${}\leftii{}$, as defined in~\eqref{action-full},
%% for uniformity).
%% \footnote{It also follows that
%%   \begin{equation*}
%%     E^{i}\vartriangleright \lambda^{n} =
%%     \q^{\half i (i - 1) - \half i n}
%%     \smash[b]{\prod_{j=0}^{i-1}} \qint{\fffrac{n}{2} - j}\;
%%     z^{i} \lambda^{n},\qquad
%%     F^{i}\vartriangleright \lambda^{n} =
%%     (-1)^i \q^{\half i (i - 1) + \half i n}
%%     \smash[b]{\prod_{j=0}^{i - 1}} \qint{\fffrac{n}{2} + j}\;
%%     \Dz^{i}\,\lambda^{n}.
%%   \end{equation*}
%% }
As we have already noted (and as is very clearly seen now), the action
restricts to $\bHDBstar$ and then pushes forward to $\HresSL2$.
There, it restricts to the subalgebra $\Czd$, and the isomorphism
\begin{equation*}
  \Czd\cong\Mat_p(\oC)
\end{equation*}
is actually that of $\UresSL2$-module
algebras~\cite{[S-U]}.%%%\enlargethispage{\baselineskip}

\subsubsection{}\label{sec:Czd-decomp}Under the above action, $\Czd$
decomposes into indecomposable $\UresSL2$ representations
as~\cite{[S-U]}
\begin{equation}\label{Czd-decomp}
  \Czd=\repP^+_1\oplus\repP^+_3\oplus\dots\oplus\repP^+_\nu,
\end{equation}
where $\nu=p-1$ if $p$ is even and $\nu=p$ if $p$ is odd, and where
$\repP^+_r$ is the projective cover of the $\UresSL2$ irreducible
representation with weight~$\q^{r-1}$ (in particular, $\repP^+_1$ is
the cover of the trivial representation; see~\cite{[FGST],[FGST2]} for
a detailed description).  The $2p$-dimensional projective module
$\repP^+_1$ in~\eqref{Czd-decomp} has the remarkable structure
\begin{equation}\label{P1}
  \xymatrix@=12pt@C=6pt{
    &&&&&\sum\limits_{i=1}^{p-1}\fffrac{1}{\qint{i}}\,z^i\,\Dz^i
    \ar^(.7){F}[dr]
    \ar_(.7){E}[dl]
    \\  
    %% *{x\Dp{p-1}\Dy^{p-1}}
    z^{p-1}
    &\kern-10pt\rightleftarrows\kern-10pt
    &
    %% *{F(x\Dp{p-1}\Dy^{p-1})}
    z^{p-2}
    &
    \kern-10pt\rightleftarrows \ldots \rightleftarrows\kern-10pt
    &z
    \ar_{F}[dr]
    &{}
    &\Dz
    \ar^{E}[dl]
    &\kern-10pt\rightleftarrows \ldots \rightleftarrows\kern-10pt
    &
    %% *{E(y\Dp{p-1}\Dy^{p-1})}
    \Dz^{p-2}
    &\kern-10pt\rightleftarrows\kern-10pt
    &
    %% *{y\Dp{p-1}\Dy^{p-1}}
    \Dz^{p-1}
    \\
    &&&&&
    1
  }
\end{equation}
where the horizontal left--right arrows denote the action of $E$ (to
the left) and $F$ (to the right) up to nonzero factors (and there are
no maps inverse to the tilted arrows).

As regards all of $\HresSL2$, its decomposition into indecomposable
$\UresSL2$ representations involves not just the ``odd'' projective
modules as in~\eqref{Czd-decomp} but actually \textit{all projective
  $\UresSL2$ modules with the multiplicity of each equal to the
  dimension of its irreducible quotient}:
\begin{equation}\label{decomp-full}
  \HresSL2=\bigoplus_{n=1}^{p}n\,\repP^+_n\oplus
  \bigoplus_{n=1}^{p}n\,\repP^-_n,
\end{equation}
where $\repP^-_r$ is the projective cover of the irreducible
representation with weight~$-\q^{r-1}$.  The multiplicities
in~\eqref{decomp-full} are identical to those in the regular
representation decomposition.\footnote{Interestingly, the sum of
  projective modules with multiplicities in the right-hand side
  of~\eqref{decomp-full} thus admits two different algebraic
  structures, one of which is actually a Hopf algebra and the other
  its module algebra.}  The sum in~\eqref{Czd-decomp} is nothing but
the $\lambda$-independent subalgebra in~$\HresSL2$.

Decomposition~\eqref{decomp-full} follows by first noting the evident
fact that the $\UresSL2$ action on $\HresSL2$ does not change the
degree in~$\lambda$, and then proceeding much as in~\cite{[S-U]}.  For
example, one of the two copies of $\repP^+_2$ involved
in~\eqref{decomp-full} is given~by
\begin{equation}\label{P2}
  \xymatrix@=12pt@C=6pt{
    &&&t_+
    \ar_(.7){E}[dl]&\kern-10pt\rightleftarrows\kern-10pt
    &t_-\ar^(.7){F}[dr]    
    \\  
    l_{p-2}
    &
    \kern-10pt\rightleftarrows \ldots \rightleftarrows\kern-10pt
    &l_1
    \ar_{F}[dr]
    &&&
    &r_1
    \ar^{E}[dl]
    &\kern-10pt\rightleftarrows \ldots \rightleftarrows\kern-10pt
    &
    r_{p-2}
    \\
    &&&
    b_+&\kern-10pt\rightleftarrows\kern-10pt&b_-
  }
\end{equation}
where
\begin{gather*}
  t_+=  \ffrac{1}{\q^2 + 1}
  \sum_{i=1}^{p - 2}\alpha_{i}\,C_{i}\,\lambda z^{i + 1}\,\Dz^{i}
  \\[-6pt]
  \intertext{with}
  \alpha_{i} = 
  \smash[t]{\sum_{j=1}^{i}\ffrac{\q^{j + \half}}{\qint{j - \half}}},
  \quad
  C_{i} = \q^{\frac{i}{2}}\smash[t]{\prod_{n=1}^{i}
    \ffrac{\qint{n - \half}}{\qint{n}}},
  \\[-6pt]
  \intertext{and}
  l_1 = \ffrac{\q^2}{\q^2 + 1}
  \smash[t]{\sum_{i=0}^{p-3}}C_{i}\,\lambda\,z^{i+2}\,\Dz^{i},
  \quad
  b_+ = \smash[t]{\sum_{i=0}^{p - 2}} C_{i}\,\lambda z^{i + 1}\,\Dz^{i}.
\end{gather*}
This construction, being a linear-in-$\lambda$ analogue of~\eqref{P1},
does not fully share its utmost simplicity, except possibly at one
point: $l_{p-2}$ in~\eqref{P2} is proportional to $\lambda z^{p-1}$;
in the other copy of $\repP^+_2$ in~\eqref{decomp-full}, linear in
$\lambda^{-1}$, $r_{p-2}$ is proportional to $\lambda^{-1}\Dz^{p-1}$.

We also note that the subspace of degree $p$ in $\lambda$ decomposes
into the sum $\repP^-_1\oplus\repP^-_3\oplus\dots\oplus\repP^-_\nu$ of
$\repP^-_{2r+1}$ modules with multiplicities $1$; because
$\lambda^{2p}=1$, there is the \textit{subalgebra}
\begin{equation*}
  \Czd+\lambda^p\Czd=
  \repP^+_1\oplus\repP^-_1\oplus\repP^+_3\oplus\repP^-_3\oplus
  \dots\oplus\repP^+_\nu\oplus\repP^-_\nu
\end{equation*}
of \textit{all} ``odd'' projective modules in~$\HresSL2$.

\subsubsection{}\label{Czd-diff}We also recall from~\cite{[S-U]} that
$\Czd$ extends to a \textit{differential} $\UresSL2$-module algebra
$\Omega\Czd$ (a quantum de~Rham complex of $\Czd$), which is the
unital algebra with the generators $z$, $\Dz$, $\dz$, $\dDz$ and the
relations (in addition to~\eqref{zp-and-dp} and~\eqref{dz-to-zd})
\begin{alignat*}{4}%%%\label{super-plane}
  \dz\,\dz&=0,&\quad \dDz\,\dDz&=0,&\quad
  \dDz\,\dz&=-\q^{-2}\dz\,\dDz,\\
  &&&&\dz\,z&=\q^{-2}z\,\dz,& \dDz\,\Dz&=\q^{2}\Dz\,\dDz,\\
  &&&&\dz\,\Dz&=\q^{2}\Dz\,\dz,& \dDz\,z&=\q^{-2}z\,\dDz.
\end{alignat*}
The differential acting as
\begin{equation*}%%%\label{the-differential}
  d(z)=\dz,\quad d(\Dz)=\dDz,\quad d(\dz)=0,\quad d(\dDz)=0
\end{equation*}
(and $d(1)=0$) 
%% The first line here immediately implies that
%% \begin{gather*}
%%   d(z^m)=\q^{1-m}\qint{m}z^{m-1}\dz,\qquad
%%   d(\Dz^n)=\q^{n-1}\qint{n}\Dz^{n-1}\dDz.
%% \end{gather*}
commutes with the $\UresSL2$ action if this is defined on $\dz$ and
$\dDz$ as
\begin{equation}\label{action-on-Omega}
  \begin{alignedat}{3}
    E\leftii\dz &= -\qint{2}z\,\dz,&\qquad
    k^2\leftii\dz &= \q^2 \dz,&\qquad
    F\leftii\dz &= 0,\\
    E\leftii\dDz &= 0,&
    k^2\leftii\dDz &= \q^{-2} \dDz,&
    F\leftii\dDz &= -\q^{2}\qint{2}\Dz\,\dDz
  \end{alignedat}
\end{equation}
and is then extended to all of $\Omega\Czd$ in accordance with the
module algebra property.
%% \begin{align*}
%%   E^{i}(z^{m}\,\dz) &= (-1)^i \q^{i m + \frac{i(i-1)}{2}}
%%   \qbin{m + i + 1}{i}\qfac{i}\, z^{m + i}\,\dz,\\
%%   F^{i}(z^{m}\,\dz) &= \q^{i(1 - m) + \frac{i(i-1)}{2}}
%%   \qbin{m}{i}\qfac{i}\,z^{m - i}\,\dz,\\
%%   E^{i}(\Dz^{m}\,\dDz) &= \q^{-i(m + 1) + \frac{i(i-1)}{2}}
%%   \qbin{m}{i}\qfac{i}\,\Dz^{m - i}\,\dDz,\\
%%   F^{i}(\Dz^{m}\,\dDz) &= (-1)^i \q^{i(m + 2) + \frac{i(i-1)}{2}}
%%   \qbin{m + i + 1}{i}\qfac{i}\,\Dz^{m + i}\,\dDz.
%% \end{align*}

%% \subsubsection{}
In fact, \textit{the entire $\HresSL2$ extends to a differential
  $\UresSL2$-module algebra}.  Let $\Omega\HresSL2$ be the algebra on
$z$, $\Dz$, $\lambda$, $\dz$, $\dDz$, and $\dlambda$ with the
relations given by those in $\Omega\Czd$ and $\HresSL2$ and the
following ones:
\begin{itemize}
\item[]$d(\lambda)=\dlambda$, \ $\dlambda\,\dlambda=0$,

\item[]$\dlambda$ commutes with $z$ and $\Dz$ and anticommutes with
  $\dz$ and $\dDz$,

\item[] $\dlambda\,\lambda= \q^{-1}\lambda\,\dlambda$ \ (whence, in
  particular, $d(\lambda^{2 p})=0$),

\item[] $\lambda$ commutes with $\dz$ and $\dDz$.
\end{itemize}
Then $\Omega\HresSL2$ endowed with the $\UresSL2$ action
\begin{equation*}
  E\leftii \dlambda = 
  \ffrac{1}{\q+1}(z\,\dlambda + \lambda\,\dz),\quad
  k^{2}\leftii \dlambda = \q^{-1} \dlambda,\quad
  F\leftii\dlambda = 
  -\ffrac{\q}{\q+1} (\Dz\,\dlambda + \lambda\,\dDz)
\end{equation*}
is a differential $\UresSL2$-module algebra.

\section{Conclusion}
We expect not only the Drinfeld double $\DD(B)$ but also the pair
$(\DD(B),\HD(B^*))$, with $\HD(B^*)$ being a $\DD(B)$-module algebra,
to play a fundamental role on the quantum group side of the
logarithmic Kazhdan--Lusztig duality.  Based on the general recipe in
Sec.~\ref{sec:thm}, the contents of Sec.~\ref{sec:sl2} must have a
counterpart for the quantum group $\mathfrak{g}_{p,p'}$ that is
Kazhdan--Lusztig-dual to the $(p,p')$ logarithmic conformal field
models~\cite{[FGST-q]}; hopefully, a ``truncation'' of the appropriate
Drinfeld double would also allow its dual version for the
corresponding Heisenberg double, yielding the pair
$(\mathfrak{g}_{p,p'},\mathfrak{h}_{p,p'})$, where
$\mathfrak{h}_{p,p'}$ is a $\mathfrak{g}_{p,p'}$-module algebra.

%% The $\UresSL2$ action on the differential module algebra
%% $\Omega\oC_{\q}[z,\Dz]\subset\HresSL2$ may also be compared to the
%% (small) quantum $s\ell(2)$ action on the de~Rham complex of the finite
%% quantum plane~\cite{[WZ]}: there, the differential is known to lift to
%% the (dual) quantum group $SL_q(2)$~\cite{[Malt],[M]} (which coacts on
%% the quantum plane).  A similar construction may also exist in our
%% case.

This paper was finished in the very inspiring atmosphere of the LCFT{}
meeting at ETH, 
\begin{wrapfigure}[7]{rH}{.21\textwidth}
  \vspace*{-.6\baselineskip}
  \includegraphics[%%%bb=68pt 142pt 870pt 880pt,
  bb=30pt 34pt 600pt 600pt,
  clip,
  width=.25\textwidth]{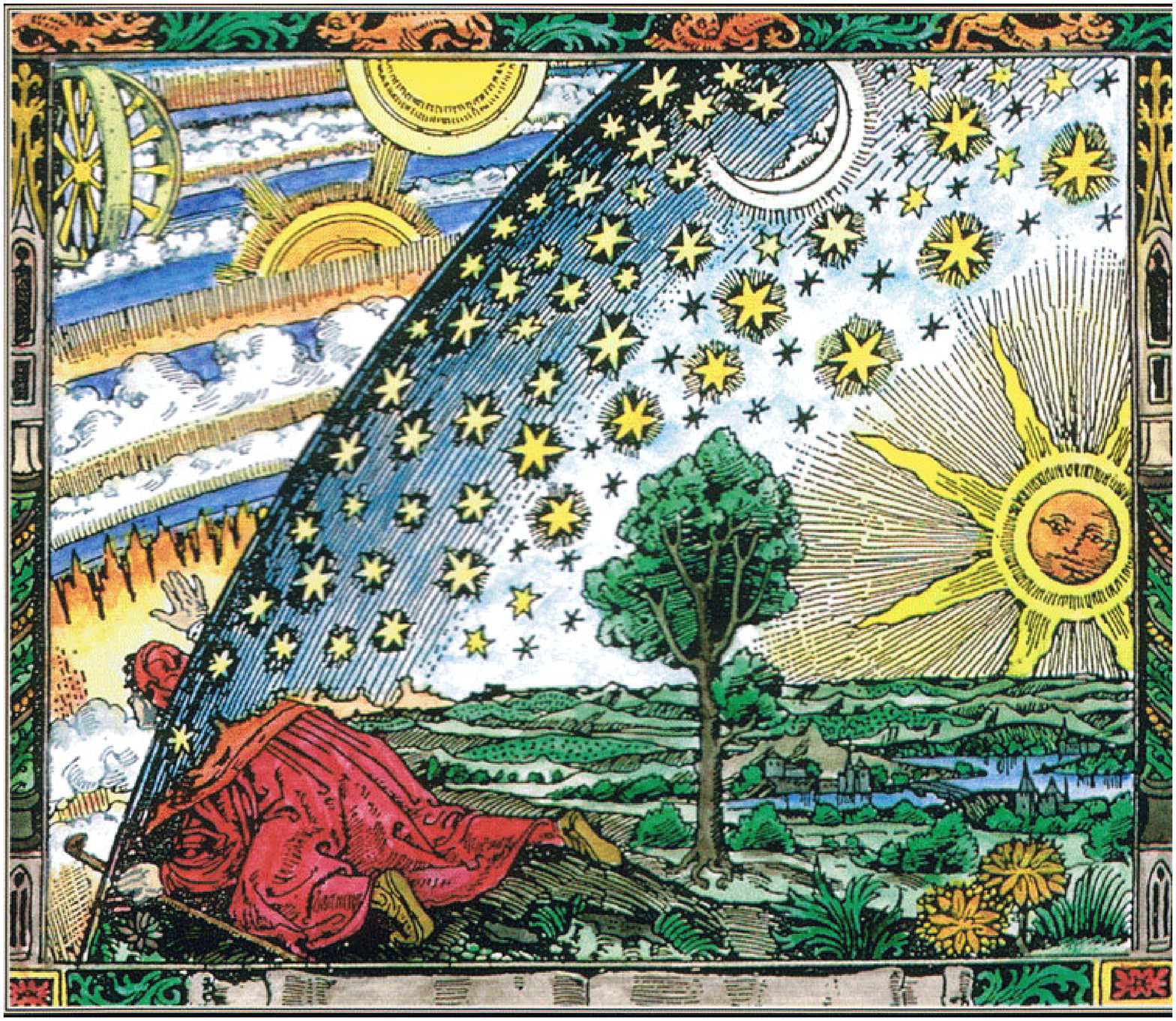}
\end{wrapfigure}
Zurich (May 2009), and it is a pleasure to thank M.~Gaberdiel for the
kind hospitality.  I~am grateful to J.~Fuchs, A.~Gainutdinov,
V.~Gurarie, P.~Mathieu, J.~Rasmussen, P.~Ruelle, I.~Runkel, and
C.~Schweigert for stimulating discussions.  Special thanks, also for
stimulation, go to G.~Mutafyan.  I thank the referee for the useful
comments.  This work was supported in part by the RFBR grant
07-01-00523, the RFBR--CNRS grant 09-01-93105, and the grant
LSS-1615.2008.2.

\appendix
\section{Drinfeld double}\label{app:D-double}
%% \subsection{}
We recall that for a Hopf algebra $B$ with bijective antipode, its
Drinfeld double $\DD(B)$ is $B^*\tensor B$ as a vector space, endowed
with the structure of a quasitriangular Hopf algebra as follows.  The
coalgebra structure is that of $B^{*\mathrm{cop}}\tensor B$, the
algebra structure is
\begin{equation}\label{in-double}
  (\mu\tensor m)(\nu\tensor n)
  = \mu(m'\leftreg \nu\rightreg S^{-1}(m'''))
  \tensor m''n
\end{equation}
for all $\mu,\nu\in B^*$ and $m,n\in B$, the antipode is given by
\begin{equation}\label{antipode-double}
  S_{_{\DD}}(\mu\tensor m)
  =(\varepsilon\tensor S(m))(\Sinv(\mu)\tensor1)
  =(S(m''')\leftreg\Sinv(\mu)\rightreg m')\tensor S(m''),
\end{equation}
and the universal $R$-matrix is
\begin{equation}\label{R-double}
  R= \sum_I(\varepsilon \tensor e_I)\tensor (e^I\tensor 1),
\end{equation}
where $\{e_I\}$ is a basis of $B$ and $\{e^I\}$ its dual basis
in~$B^*$.

\section{LCFT{} motivation}\label{app:lcftrem}
Referring to the $(p,1)$ logarithmic conformal models, we here
emphasize several features that may be captured on the algebraic side
by $\Czd$\,---\,the ``noncommutative part'' of $\HresSL2$\,---\,and
its de~Rham complex $\Omega\Czd$ (Secs.~\bref{sec:Czd-decomp}
and~\bref{Czd-diff}).  The basic observation is that the fields of the
simplest, $(p=2,1)$ logarithmic model arrange under the action of $E$
and $F$ into the same diagrams as certain elements of $\Omega\Czd$ at
$p=2$ (i.e., $\q=\sqrt{-1}$).  Guided by the quantum group symmetry,
we then expect that for $p>2$, the corresponding diagrams for
$\Omega\Czd$ capture some structures to be observed in the $(p,1)$
models when these are described in terms of fields manifestly
covariant under $\UresSL2$.  The quantum group analysis may thus help
find this (so far hypothetical) description of logarithmic conformal
models.

We proceed from the free-fermion description of the $(p=2,1)$
logarithmic conformal field model.  The starting point is the usual
system of two free fermion fields $\xi(u)$ and $\eta(u)$ with the
respective conformal weights $0$ and $1$, whose OPE is
\begin{equation*}
  \xi(u)\,\eta(v)=\ffrac{1}{u-v},\qquad u,v\in\oC.
\end{equation*}
Virasoro generators with central charge $c=-2$ are the modes of the
energy--momentum tensor
\begin{equation*}
  T(u)=-\eta(u)\dd\xi(u)=\sum_{n\in\oZ}L_n z^{-n-2},
\end{equation*}
where $\dd=\partial/\partial u$ and the normal-ordered product is
understood in the right-hand side.  It follows that two
screenings\,---\,operators commuting with this Virasoro
algebra\,---\,are given by
\begin{align}\label{screening}
  E&=\oint\eta=\eta_0\\
  \intertext{(the ``short'' screening, which squares to zero for
    $p=2$) and}
\label{screening-long}
  \mathsf{f}&=\oint\dd\xi\,\xi
\end{align}
(the ``long'' screening).  The relevant complex of (Feigin--Fuchs)
Virasoro modules is
\begin{equation}\label{FF}
  \xymatrix@=28pt@C=46pt{
    &&\Socle{1}="1s"&\Quotient{\xi}="m1ms"\\
    &\Socle{\eta}="1ms"&\Quotient{\quad\eta\xi}="1s2"&\Socle{}="m1ms2"
    &\Quotient{\dd\xi\,\xi}="m1s"\\
    \Socle{\dd\eta\eta}="3s"&\Quotient{}="1ms2"&\Socle{}="1s3"&\Quotient{}="m1ms3"&\Socle{\qquad\quad\dd^2\xi\dd\,\xi}="m1s2"&\Quotient{\quad\dd^2\xi\dd\,\xi\,\xi}="m3ms"\\
    \Quotient{}="3s2"&\Socle{}="1ms3"&\Quotient{}="1s4"&\Socle{}="m1ms4"&\Quotient{}="m1s3"&\Socle{}="m3ms2"&\Quotient{}="m3s"
%%     &&&&&&&&&
    %%%%%%%%%%% SCREENINGS
    \ar"1ms2";"3s"_E
    \ar"1s2";"1ms"_E
    \ar"m1ms";"1s"_E
    \ar"m1s";"m1ms2"_E
    \ar"m3ms";"m1s2"^E
    \ar"m3s";"m3ms2"_E
    \ar"1s4";"1ms3"_E
    \ar"m1ms3";"1s3"_E
    \ar"m1s3";"m1ms4"_E
    %%%%%%%%%%% SINGULAR AND COSINGULAR
    \ar"1ms2";"1ms"
    \ar"1ms2";"1ms3"
    \ar"3s2";"3s"
    \ar"1s2";"1s"
    \ar"1s2";"1s3"
    \ar"1s4";"1s3"
    \ar"m1ms";"m1ms2"
    \ar"m1ms3";"m1ms2"
    \ar"m1ms3";"m1ms4"
    \ar"m1s";"m1s2"
    \ar"m1s3";"m1s2"
    \ar"m3ms";"m3ms2"
    \ar"3s2";"3s2"+<-30pt,0pt>
  }
\end{equation}

\smallskip

\noindent
where vertical arrows (directed towards \textit{sub}modules) indicate
embedding of subquotients in Feigin--Fuchs modules.  The picture
continues to the left and to the right (and downward) indefinitely.
The weight-$2$ fields $\dd\eta\,\eta$ and $\dd^2\xi\,\dd\xi$ are the
triplet algebra generators.

The algebra of fields is then extended by a field $\ddinv\eta(u)$ such
that
\begin{equation*}
  \xymatrix@=12pt@R=16pt{
    \ddinv\eta(u)\ar^{L_{-1}}[1,0]\\
    \eta(u)
  }
\end{equation*}
It is $\delta(u)=\ddinv\eta(u)$ and $\xi(u)$ that are in fact the
symplectic fermions~\cite{[K-sy]}.  These \textit{weight-$0$} fields
generate two standard first-order systems: our starting
$(\eta(u),\xi(u))$ and $(\delta(u),$\linebreak[0]$\dd\xi(u))$
(cf.~\cite{[FGST2]}).

The fields $(\delta(u),\xi(u))$ allow constructing a logarithmic
partner $\Lambda(u)=\delta(u)\,\xi(u)$ of the identity operator;
diagram~\eqref{FF} then extends such that the top level (split
vertically for visual clarity) becomes
\begin{equation}\label{split}
  \xymatrix@=12pt@R=-8pt{%
    &&&&\delta(u)\,\xi(u)\ar[1,4]^F\ar[1,-4]_E&&&&&\\
    \delta(u)\ar[1,4]_F&&&&&&&&\xi(u)\ar[1,-4]^E&\\
    &&&&1&&&&&
  }
\end{equation}

There are also two characteristic diagrams of \textit{weight-$1$}
fields.  We recall that if the fermions are bosonized through a free
bosonic field,
\begin{equation*}
  \xi(u)=e^{\varphi(u)},\quad  \eta(u)=e^{-\varphi(u)},
  \quad
  \eta(u)\xi(u)=-\dd\varphi(u),
\end{equation*}
then the long-screening \textit{current} (the ``integrand''
in~\eqref{screening-long}) is $e^{2\varphi}$ (which is a weight-$1$
field), and we have
\begin{equation}\label{scr-right}
  \xymatrix@=12pt{%
    \delta(u)\dd\xi(u)\ar[1,1]^F&&e^{2\varphi(u)}\ar[1,-1]_E\\
    &\dd\xi(u)
  }
\end{equation}
Similarly, there is an alternative bosonization through the scalar
field introduced as $\dd\phi(u)=\delta(u)\dd\xi(u)$.  This gives the
diagram
\begin{equation}\label{scr-left}
  \xymatrix@=12pt{%
    e^{2\phi(u)}\ar[1,1]^F&&\eta(u)\xi(u)\ar[1,-1]_E\\
    &\eta(u)
  }
\end{equation}
(once again, $\eta(u)=\dd\delta(u)$, which makes the two diagrams
symmetric to each other).

The $(p=2,1)$ logarithmic model corresponds to $\q=\sqrt{-1}$
in~\eqref{the-q}.  Relations~\eqref{zp-and-dp} and~\eqref{dz-to-zd} in
$\Czd$ are then indeed those mimicking free fermions:
\begin{equation*}
  z^2=0,\quad \Dz^2=0,\quad
  \Dz z + z\Dz = 2i.
\end{equation*}
Based on the $\UresSL2$ symmetry, we conjecture that for general $p$,
$\Czd$ similarly allows expressing the relations among the
$(p,1)$-model fields ``with an explicit quantum-group index.''
Following~\cite{[FGST2]}, we call such fields \textit{parafermions}
(the term is somewhat overloaded by different meanings; its usage for
fields transforming under a quantum group action goes back
to~\cite{[Sm]}).

On the quantum-group side, clearly, \eqref{P1} is the general-$p$
counterpart of~\eqref{split} under the correspondence
\begin{equation*}
  z^i\leftrightarrow\xi^{(i)}(u),\qquad
  \Dz^j\leftrightarrow\delta^{(j)}(u)
\end{equation*}
for the $(p-1)$-component ``parafermion'' fields $\delta^{(i)}(u)$ and
$\xi^{(i)}(u)$ generalizing the symplectic fermions
$(\delta(u),\xi(u))$.  The constituents of~\eqref{P1} satisfy
commutation relations generalizing the fermionic ones that occur for
$p=2$: for general $p$, we have
\begin{gather*}
  \Dz^j\,z^i=\sum_{\ell\geq0}\q^{-(2 j - \ell) i + \ell j - \frac{\ell(\ell-1)}{2}}
  \qbin{j}{\ell} \qbin{i}{\ell} \qfac{\ell} \left(\q - \q^{-1}\right)^\ell
  z^{i-\ell}\Dz^{j-\ell}.
\end{gather*}

Moreover, the counterparts of~\eqref{scr-right} and~\eqref{scr-left}
for general $p$ are the diagrams that are easily established
using~\eqref{action-on-Omega}, essentially by applying the
differential to~\eqref{P1}, with the resulting modules extended by the
``cohomology corners'' $z^{p-1}\dz$ and~$\Dz^{p-1}\dDz$:
\begin{gather*}
  \xymatrix@=12pt@C=5pt{
    \displaystyle\smash[t]{\sum_{i=1}^{p - 1}}
    \fffrac{1}{\qint{i}}z^{i}\,d(\Dz^{i})
    \ar^(.6)F[dr]
    &&&&&&\Dz^{p-1}\,\dDz\ar_E[dl]
    \\
    &\dDz
    &\kern-10pt\rightleftarrows\kern-10pt
    &\Dz\,\dDz
    &\kern-6pt\rightleftarrows \ldots \rightleftarrows\kern-6pt
    &\Dz^{p-2}\,\dDz
  }
  \\
  \intertext{and}
%%%\label{two-modules}
  \xymatrix@=12pt@C=5pt{
    z^{p-1}\,\dz
    \ar^F[dr]&&&&&&
    %% \displaystyle\sum_{i=1}^{p - 1}\q^{1 - i} z^{i - 1}\dz\Dz^{i}    
      \displaystyle\smash[t]{\sum_{i=1}^{p - 1}}
      \fffrac{1}{\qint{i}}d(z^{i})\,\Dz^{i}
    \ar_(.75)E[dl]
    \\
    &z^{p-2}\,\dz
    &\kern-6pt\rightleftarrows \ldots \rightleftarrows\kern-6pt
    &z\,\dz
    &\kern-10pt\rightleftarrows\kern-10pt
    &\dz&
  }
\end{gather*}
(as before, horizontal left--right arrows represent the action of $E$
and $F$ up to nonzero factors and tilted arrows have no inverse maps).
The cohomology corners are Hopf-algebra counterparts of the screening
currents in the two bosonizations, and the bottom elements, of the
differentials $\dd\delta^{(1)}(u)$, \dots, $\dd\delta^{(p-1)}(u)$ and
$\dd\xi^{(p-1)}(u)$, \dots, $\dd\xi^{(1)}(u)$.

\parindent=0pt

\end{document}